\documentclass[11pt]{amsart}%
\usepackage{amsmath}
\usepackage{graphicx}
\usepackage{amsfonts}
\usepackage{amssymb}%
\newtheorem{theorem}{Theorem}
\theoremstyle{plain}

\newtheorem{corollary}[theorem]{Corollary}

\newtheorem{example}[theorem]{Example}

\newtheorem{lemma}[theorem]{Lemma}

\newtheorem{proposition}[theorem]{Proposition}

\begin{document}
\title{Extension of Functions with Small Oscillation}

\begin{abstract}
A classical theorem of Kuratowski says that every Baire one function on a
$G_{\delta}$ subspace of a Polish ($=$ separable completely metrizable) space
$X$ can be extended to a Baire one function on $X$. Kechris and Louveau
introduced a finer gradation of Baire one functions into small Baire classes.
A Baire one function $f$ is assigned into a class in this heirarchy depending
on its oscillation index $\beta(f)$. We prove a refinement of Kuratowski's
theorem: if $Y$ is a subspace of a metric space $X$ and $f$ is a real-valued
function on $Y$ such that $\beta_{Y}(f)<\omega^{\alpha}$, $\alpha<\omega_{1}$,
then $f$ has an extension $F$ onto $X$ so that $\beta_{X}(F)\leq\omega
^{\alpha}$. We also show that if $f$ is a continuous real valued function on
$Y,$ then $f$ has an extension $F$ onto $X$ so that $\beta_{X}\left(
F\right)  \leq3.$ An example is constructed to show that this result is optimal.

\end{abstract}
\author{Denny H. Leung}
\address{Department of Mathematics, National University of Singapore, 2 Science Drive
2, Singapore 117543.}
\email{matlhh@nus.edu.sg}
\author{Wee-Kee Tang}
\address{Mathematics and Mathematics Education, National Institute of Education \\
Nanyang Technological University, 1 Nanyang Walk, Singapore 637616.}
\email{wktang@nie.edu.sg}
\maketitle

Let $X$ be a topological space. A real-valued function on $X$ belongs to Baire
class one if it is the pointwise limit of a sequence of continuous functions.
If $X$ is a Polish ($=$ separable completely metrizable) space, then a
classical theorem of Kuratowski \cite{K} states that every Baire one function
on a $G_{\delta}$ subspace of $X$ can be extended to a Baire one function on
$X$. In \cite{KL}, Kechris and Louveau introduced a finer gradation of Baire
one functions into small Baire classes using the oscillation index $\beta$,
whose definition we now recall.

Let $X$ be a topological space and let ${\mathcal{C}}$ denote the collection
of all closed subsets of $X$. A \emph{derivation} on ${\mathcal{C}}$ is a map
${\mathcal{D}} : {\mathcal{C}} \to{\mathcal{C}}$ such that ${\mathcal{D}}(P)
\subseteq P$ for all $P \in{\mathcal{C}}$. The oscillation index $\beta$ is
associated with a family of derivations. Let $\varepsilon> 0$ and a function
$f : X \to{\mathbb{R}}$ be given. For any $P \in{\mathcal{C}}$, let
${\mathcal{D}}^{1}(f, \varepsilon, P)$ be the set of all $x \in P$ such that
for any neighborhood $U$ of $x$, there exist $x_{1}, x_{2} \in P \cap U$ such
that $|f(x_{1}) - f(x_{2})| \geq\varepsilon$. The derivation ${\mathcal{D}%
}^{1}(f,\varepsilon,\cdot)$ may be iterated in the usual manner. For all
$\alpha< \omega_{1}$, let
\[
{\mathcal{D}}^{\alpha+1}(f, \varepsilon, P) = {\mathcal{D}}^{1}(f,
\varepsilon, {\mathcal{D}}^{\alpha}(f, \varepsilon, P)).
\]
If $\alpha$ is a countable limit ordinal, set
\[
{\mathcal{D}}^{\alpha}(f, \varepsilon, P) = \cap_{\gamma<\alpha}{\mathcal{D}%
}^{\gamma}(f, \varepsilon, P).
\]
If ${\mathcal{D}}^{\alpha}(f,\varepsilon,P) \neq\emptyset$ for all $\alpha<
\omega_{1}$, let $\beta_{X}(f,\varepsilon) = \omega_{1}$. Otherwise, let
$\beta_{X}(f,\varepsilon)$ be the smallest countable ordinal $\alpha$ such
that ${\mathcal{D}}^{\alpha}(f,\varepsilon,P) = \emptyset$. The
\emph{oscillation index} of $f$ is $\beta_{X}(f) = \sup_{\varepsilon> 0}%
\beta_{X}(f,\varepsilon)$.

The main result of \S 1 is that if $Y$ is a subspace of a metric space $X$ and
$f:Y\rightarrow{\mathbb{R}}$ satisfies $\beta_{Y}(f)<\omega^{\alpha}$ for some
$\alpha<\omega_{1}$, then $f$ can be extended to a function $F$ on $X$ with
$\beta_{X}(F)\leq\omega^{\alpha}$. It follows readily from the Baire
Characterization Theorem \cite[10.15]{BBT} that a real-valued function on a
Polish space is Baire one if and only if its oscillation index is countable.
(See, e.g. \cite{KL}.) Also, a theorem of Alexandroff says that a $G_{\delta}$
subspace of a Polish space is Polish \cite[10.18]{BBT}. Hence our result
refines Kuratowski's theorem in terms of the oscillation index. Let us mention
that if $X$ is a metric space, then every real-valued function with countable
oscillation index on a closed subspace of $X$ may be extended onto $X$ with
preservation of the index \cite[Theorem 3.6]{LT}. (Note that the proof of
\cite[Theorem 3.6]{LT} does not require the compactness of the ambient space.)
More recent results on the extension of Baire one functions on general
topological spaces are found in \cite{KS}.

It is well known that if a function is continuous on a \emph{closed} subspace
of a \ metric space, then there exists a continuous extension to the whole
space. \S 2 is devoted to the study of extensions of continuous functions from
an \emph{arbitrary} subspace of a metric space. It is shown that if $f$ is a
continuous function on a subspace $Y$ of a metric space $X,$ then $f$ has an
extension $F$ on $X$ with $\beta_{X}\left(  F\right)  \leq3.$ An example is
given to show that the result is optimal. The criteria for continuous
extension on dense subspaces had been studied by several authors. (See, e.g.,
\cite{B}, \cite{H}.)

\section{Functions of Small Oscillation}

Given a real-valued function defined on a set $S,$ let $\left\Vert
f\right\Vert _{S}=\sup_{s\in S}\left\vert f\left(  s\right)  \right\vert .$
For any topological space $X,$ the support $\operatorname*{supp}f$ of a
function $f:X\rightarrow\mathbb{R}$ is the closed set $\overline{\left\{  x\in
X:f\left(  x\right)  \neq0\right\}  }.$ A family $\left\{  \varphi_{\alpha
}:\alpha\in\mathcal{A}\right\}  $ of nonnegative real-valued functions on $X$
is called a \emph{partition of unity on }$X$ if

\begin{enumerate}
\item The support of the $\varphi_{\alpha}$'s form a locally finite closed
covering of $X,$

\item $\sum_{\alpha\in\mathcal{A}}\varphi_{\alpha}\left(  x\right)  =1$ for
all $x\in X.$
\end{enumerate}

\noindent If $\left\{  U_{\beta}:\beta\in\mathcal{B}\right\}  $ is an open
covering of $X,$ we say that a partition of unity $\left\{  \varphi_{\beta
}:\beta\in\mathcal{B}\right\}  $ on $X$ is subordinated to $\left\{  U_{\beta
}:\beta\in\mathcal{B}\right\}  \,$if the support of each $\varphi_{\beta}$
lies in the corresponding $U_{\beta}.$ It is well known that if $X$ is
paracompact (in particular, if $X$ is a metric space \cite[Theorem IX 5.3]%
{D}), then for each open covering $\left\{  U_{\beta}:\beta\in\mathcal{B}%
\right\}  $ of $X$ there is a partition of unity on $X$ subordinated to
$\left\{  U_{\beta}:\beta\in\mathcal{B}\right\}  .$ (See, for example,
\cite[Theorem VIII 4.2]{D}.)

\begin{proposition}
\label{P1}Let $X$ be a metric space and $Y$ be a subspace of $X.$ Suppose that
$f:Y\rightarrow\mathbb{R}$ is a function such that $\beta_{Y}\left(
f,\varepsilon\right)  \leq\alpha$ for some $\varepsilon>0,$ $\alpha<\omega
_{1}.$ Then there exists a function $\tilde{f}:X\rightarrow\mathbb{R}$ such
that $\beta_{X}\left(  \tilde{f}\right)  \leq\alpha+1,$ $\left\Vert \tilde
{f}\right\Vert _{X}\leq\left\Vert f\right\Vert _{Y}$ and $\left\Vert
f-\tilde{f}\right\Vert _{Y}\leq\varepsilon.$
\end{proposition}

In the following, denote $\mathcal{D}^{\beta}\left(  f,\varepsilon,Y\right)  $
by $Y^{\beta}$ for all $\beta<\omega_{1}.$ Proposition \ref{P1} is proved by
working on each of the pieces $Y^{\beta}\smallsetminus Y^{\beta+1}$,
$\beta<\alpha$, and gluing together the results.

\begin{lemma}
\label{L3}For all $0\leq\beta<\alpha,$ there exist an open set $Z_{\beta}$ in
$X$ such that $Y^{\beta}\smallsetminus Y^{\beta+1}\subseteq Z_{\beta}%
\subseteq\left(  Y^{\beta+1}\right)  ^{c},$ and a continuous function
$f_{\beta}:Z_{\beta}\rightarrow\mathbb{R}$ such that $\left\Vert f-f_{\beta
}\right\Vert _{Y^{\beta}\smallsetminus Y^{\beta+1}}\leq\varepsilon$ and
$\left\Vert f_{\beta}\right\Vert _{Z_{\beta}}\leq\left\Vert f\right\Vert
_{Y}.$
\end{lemma}

\begin{proof}
If $0\leq\beta<\alpha$ and $y\in Y^{\beta}\smallsetminus Y^{\beta+1},$ there
exists a set $U_{y}$ that is an open neighborhood of $y$ in $X$ so that
$U_{y}$ is disjoint from $Y^{\beta+1}$ and that $f\left(  U_{y}\cap Y^{\beta
}\right)  \subseteq\left(  f\left(  y\right)  -\varepsilon,f\left(  y\right)
+\varepsilon\right)  .$ Let
\[
Z_{\beta}=%
{\displaystyle\bigcup\limits_{y\in Y^{\beta}\smallsetminus Y^{\beta+1}}}
U_{y}.
\]
Each $Z_{\beta}$ is open in $X$. Clearly, $Y^{\beta}\smallsetminus Y^{\beta
+1}\subseteq Z_{\beta}\subseteq\left(  Y^{\beta+1}\right)  ^{c}$. There exists
a partition of unity $\left(  \varphi_{y}\right)  _{y\in Y^{\beta
}\smallsetminus Y^{\beta+1}}$ on $Z_{\beta}$ subordinated to the open covering
$\mathcal{U}=\left\{  U_{y}:y\in Y^{\beta}\smallsetminus Y^{\beta+1}\right\}
$. Define $f_{\beta}:Z_{\beta}\rightarrow\mathbb{R}$ by
\[
f_{\beta}\left(  z\right)  =\sum_{y\in Y^{\beta}\smallsetminus Y^{\beta+1}%
}f\left(  y\right)  \varphi_{y}\left(  z\right)  .
\]
Then $f_{\beta}$ is well-defined, continuous and $\left\Vert f_{\beta
}\right\Vert _{Z_{\beta}}\leq\left\Vert f\right\Vert _{Y}.$ If $x\in Y^{\beta
}\smallsetminus Y^{\beta+1},$ set $V_{x}=\left\{  y\in Y^{\beta}\smallsetminus
Y^{\beta+1}:\varphi_{y}\left(  x\right)  \neq0\right\}  .$ Then $\sum_{y\in
V_{x}}\varphi_{y}\left(  x\right)  =1$. If $y\in V_{x},$ then $x\in U_{y};$
thus $\left\vert f\left(  x\right)  -f\left(  y\right)  \right\vert
<\varepsilon$. Hence%
\begin{align*}
\left\vert f\left(  x\right)  -f_{\beta}\left(  x\right)  \right\vert  &
=\left\vert \sum_{y\in V_{x}}\left(  f\left(  x\right)  -f\left(  y\right)
\right)  \varphi_{y}\left(  x\right)  \right\vert \\
&  \leq\sum_{y\in V_{x}}\left\vert f\left(  x\right)  -f\left(  y\right)
\right\vert \varphi_{y}\left(  x\right)  \leq\varepsilon.
\end{align*}
Therefore, $\left\Vert f-f_{\beta}\right\Vert _{Y^{\beta}\smallsetminus
Y^{\beta+1}}\leq\varepsilon$, as required.
\end{proof}

\begin{proof}
[Proof of Proposition \ref{P1}]Define a function $\tilde{f}:X\rightarrow
\mathbb{R}$ by%
\[
\tilde{f}\left(  x\right)  =\left\{
\begin{array}
[c]{ccc}%
f_{\beta}\left(  x\right)  & \text{if} & x\in Z_{\beta}\smallsetminus
\cup_{\gamma<\beta}Z_{\gamma},\text{ }\beta<\alpha,\\
0 & \text{if} & x\notin\cup_{\gamma<\alpha}Z_{\gamma}.
\end{array}
\right.
\]
Clearly, $\left\Vert \tilde{f}\right\Vert _{X}=\sup_{\beta<\alpha}\left\Vert
f_{\beta}\right\Vert _{Z_{\beta}}\leq\left\Vert f\right\Vert _{Y}.$ If $x\in
Y,$ then $x\in Y^{\beta}\smallsetminus Y^{\beta+1}$ for some $\beta<\alpha.$
In particular, $x\in Z_{\beta}\smallsetminus\cup_{\gamma<\beta}Z_{\gamma}$.
Hence $\left\vert f\left(  x\right)  -\tilde{f}\left(  x\right)  \right\vert
=\left\vert f\left(  x\right)  -f_{\beta}\left(  x\right)  \right\vert
\leq\left\Vert f-f_{\beta}\right\Vert _{Y^{\beta}\smallsetminus Y^{\beta+1}%
}\leq\varepsilon$ according to Lemma \ref{L3}. Since this is true for all
$x\in Y,$ $\left\Vert f-\tilde{f}\right\Vert _{Y}\leq\varepsilon.$

It remains to show that $\beta_{X}\left(  \tilde{f}\right)  \leq\alpha+1.$ To
this end, we claim that $\mathcal{D}^{\beta}\left(  \tilde{f},\delta,X\right)
\cap Z_{\gamma}=\emptyset$ for all $\delta>0,$ $\gamma<\beta\leq\alpha.$ We
prove the claim by induction$.$ Let $\delta>0.$ Since $f_{0}$ is continuous on
the open set $Z_{0},$ we have $\mathcal{D}^{1}\left(  \tilde{f},\delta
,X\right)  \cap Z_{0}=\emptyset.$ Suppose that the claim holds for all
ordinals less than $\beta.$ By the inductive hypothesis, $\mathcal{D}^{\xi
}\left(  \tilde{f},\delta,X\right)  \cap\left(  \cup_{\gamma<\xi}Z_{\gamma
}\right)  =\emptyset$ for all $\xi<\beta.$ Therefore,%

\[
\mathcal{D}^{\xi}\left(  \tilde{f},\delta,X\right)  \cap\left[  Z_{\xi
}\smallsetminus\left(  \cup_{\gamma<\xi}Z_{\gamma}\right)  \right]
=\mathcal{D}^{\xi}\left(  \tilde{f},\delta,X\right)  \cap Z_{\xi}.
\]
Now $\tilde{f}=f_{\xi}$ is continuous on this set, which is open in
$\mathcal{D}^{\xi}\left(  \tilde{f},\delta,X\right)  .$ Therefore
$\mathcal{D}^{\xi+1}\left(  \tilde{f},\delta,X\right)  \cap Z_{\xi}%
=\emptyset.$ Also since $\mathcal{D}^{\beta}\left(  \tilde{f},\delta,X\right)
\subseteq\mathcal{D}^{\gamma+1}\left(  \tilde{f},\delta,X\right)  $ for all
$\gamma<\beta,$
\[
\mathcal{D}^{\beta}\left(  \tilde{f},\delta,X\right)  \cap Z_{\gamma
}=\emptyset
\]
for all $\gamma<\beta.$ This proves the claim. It follows from the claim that%
\[
\mathcal{D}^{\alpha}\left(  \tilde{f},\delta,X\right)  \subseteq\left(
\cup_{\gamma<\alpha}Z_{\gamma}\right)  ^{c}%
\]
for any $\delta>0.$ Since $\tilde{f}=0$ on the latter set, $\mathcal{D}%
^{\alpha+1}\left(  \tilde{f},\delta,X\right)  =\emptyset.$
\end{proof}

In order to iterate Proposition \ref{P1} to obtain an extension of $f,$ we
need the following result.

\begin{proposition}
\label{P2}Let $Y$ be a subspace of a metric space $X.$ If $\beta_{Y}\left(
f\right)  <\omega^{\xi}$ and $\beta_{Y}\left(  g\right)  <\omega^{\xi},$ then
$\beta_{Y}\left(  f+g\right)  <\omega^{\xi}.$
\end{proposition}

Proposition \ref{P2} is proved by the method used in \cite[Lemma 5]{KL}. This
requires a slight modification in the derivation $\mathcal{D}$ associated with
the index $\beta.$

Given a real valued function $f:Y\rightarrow\mathbb{R}$, $\varepsilon>0,$ and
a closed subset $P$ of $Y,$ define $G\left(  f,\varepsilon,P\right)  $ to be
the set of all $y\in P$ such that for every neighborhood $U$ of $y,$ there
exists $y^{\prime}\in P\cap U$ such that $\left\vert f\left(  y\right)
-f\left(  y^{\prime}\right)  \right\vert \geq\varepsilon.$ Let%

\[
\mathcal{G}^{1}\left(  f,\varepsilon,P\right)  =\overline{G\left(
f,\varepsilon,P\right)  },
\]
where the closure is taken in $Y.$ This defines a derivation $\mathcal{G}$ on
the closed subsets of $Y$ which may be iterated in the usual manner. If
$\alpha<\omega_{1},$ let%
\[
\mathcal{G}^{\alpha+1}\left(  f,\varepsilon,P\right)  =\mathcal{G}^{1}\left(
f,\varepsilon,\mathcal{G}^{\alpha}\left(  f,\varepsilon,P\right)  \right)  .
\]
If $\alpha<\omega_{1}$ is a limit ordinal, let%
\[
\mathcal{G}^{\alpha}\left(  f,\varepsilon,P\right)  =%
{\displaystyle\bigcap\limits_{\alpha^{\prime}<\alpha}}
\mathcal{G}^{\alpha^{\prime}}\left(  f,\varepsilon,P\right)  .
\]
Clearly, the derivation $\mathcal{G}$\ is closely related to $\mathcal{D}.$
The precise relationship between $\mathcal{D}$ and $\mathcal{G}$ is given in
part (c) of the next lemma.

\begin{lemma}
\label{L1}If $f$ and $g$ are real-valued functions on $Y$, $\varepsilon>0,$
and $P,Q$ are closed subsets of $Y,$ then

(a) $\mathcal{G}^{1}\left(  f+g,\varepsilon,P\right)  \subseteq\mathcal{G}%
^{1}\left(  f,\varepsilon/2,P\right)  \cup\mathcal{G}^{1}\left(
g,\varepsilon/2,P\right)  ,$

(b) $\mathcal{G}^{1}\left(  f,\varepsilon,P\cup Q\right)  \subseteq
\mathcal{G}^{1}\left(  f,\varepsilon,P\right)  \cup\mathcal{G}^{1}\left(
f,\varepsilon,Q\right)  ,$

(c) $\mathcal{D}^{1}\left(  f,2\varepsilon,P\right)  \subseteq\mathcal{G}%
^{1}\left(  f,\varepsilon,P\right)  \subseteq\mathcal{D}^{1}\left(
f,\varepsilon,P\right)  .$
\end{lemma}

We leave the simple proofs to the reader. Note that it follows from part (c)
that for all $\alpha<\omega_{1}$,
\begin{equation}
\mathcal{D}^{\alpha}\left(  f,2\varepsilon,P\right)  \subseteq\mathcal{G}%
^{\alpha}\left(  f,\varepsilon,P\right)  \subseteq\mathcal{D}^{\alpha}\left(
f,\varepsilon,P\right)  . \tag{d}\label{3'}%
\end{equation}

\begin{proof}
[Proof of Proposition \ref{P2}]Parts (a) and (b) of Lemma \ref{L1} correspond
to (*) and (**) in \cite[Lemma 5]{KL} respectively. From the proof of that
result, we obtain for all $n\in\mathbb{N}$ and $\zeta<\omega_{1},$%
\begin{equation}
\mathcal{G}^{\omega^{\zeta}\cdot2n}\left(  f+g,\varepsilon,Y\right)
\subseteq\mathcal{G}^{\omega^{\zeta}\cdot n}\left(  f,\varepsilon/2,Y\right)
\cup\mathcal{G}^{\omega^{\zeta}\cdot n}\left(  g,\varepsilon/2,Y\right)  .
\label{E1}%
\end{equation}
Since $\beta_{Y}\left(  f\right)  <\omega^{\xi}$ and $\beta_{Y}\left(
g\right)  <\omega^{\xi},$ there exist $\zeta<\xi$ and $n\in\mathbb{N}$ such
that $\beta_{Y}\left(  f\right)  <\omega^{\zeta}\cdot n$ and $\beta_{Y}\left(
g\right)  <\omega^{\zeta}\cdot n$. By (\ref{3'}), for any $\varepsilon>0,$%
\[
\mathcal{G}^{\omega^{\zeta}\cdot n}\left(  f,\varepsilon/2,Y\right)
=\mathcal{G}^{\omega^{\zeta}\cdot n}\left(  g,\varepsilon/2,Y\right)
=\emptyset.
\]
By (\ref{3'}) and (\ref{E1}) ,
\[
\mathcal{D}^{\omega^{\zeta}\cdot2n}\left(  f+g,2\varepsilon,Y\right)
=\emptyset.
\]
Since this is true for all $\varepsilon>0,$ we have
\[
\beta_{Y}\left(  f+g\right)  \leq\omega^{\zeta}\cdot2n<\omega^{\xi}.
\]

\end{proof}

\begin{theorem}
Let $X$ be a metric space and let $Y$ be an arbitrary subspace of $X.$ If
$f:Y\rightarrow\mathbb{R}$ satisfies $\beta_{Y}\left(  f\right)
<\omega^{\alpha}$ for some $\alpha<\omega_{1},$ then there exists
$F:X\rightarrow\mathbb{R}$ with $\beta_{X}\left(  F\right)  \leq\omega
^{\alpha}$ and $F_{|Y}=f.$
\end{theorem}

\begin{proof}
Applying Proposition \ref{P1} to $f:Y\rightarrow\mathbb{R}$ with
$\varepsilon=\frac{1}{2},$ we obtain $g_{1}:X\rightarrow\mathbb{R}$, with
$\left\|  f-g_{1}\right\|  _{Y}\leq\frac{1}{2},$ and $\beta_{X}\left(
g_{1}\right)  <\omega^{\alpha}$. By Proposition \ref{P2}, $\beta_{Y}\left(
f-g_{1}\right)  <\omega^{\alpha}.$ Now apply Proposition \ref{P1} to $\left(
f-g_{1}\right)  _{|Y}$ with $\varepsilon=\frac{1}{2^{2}}.$ We obtain
$g_{2}:X\rightarrow\mathbb{R}$, with $\left\|  g_{2}\right\|  _{X}\leq\left\|
f-g_{1}\right\|  _{Y}\leq\frac{1}{2},$ $\left\|  f-g_{1}-g_{2}\right\|
_{Y}\leq\frac{1}{2^{2}},$ and $\beta_{X}\left(  g_{2}\right)  <\omega^{\alpha
}.$ Continuing in this way, we obtain a sequence $\left(  g_{n}\right)  $ of
real-valued functions on $X$ such that for all $n\in\mathbb{N},$

(i) $\left\Vert g_{n+1}\right\Vert _{X}\leq\left\Vert f-\sum_{i=1}^{n}%
g_{i}\right\Vert _{Y}\leq\frac{1}{2^{n}},$

(ii) $\beta_{X}\left(  g_{n}\right)  <\omega^{\alpha}.$

\noindent Let $F=\sum_{n=1}^{\infty}g_{n}.$ Note that the series converges
uniformly on $X$ and $g_{|Y}=f$ by (i). Finally, suppose that $\varepsilon>0.$
Choose $N$ such that $\sum_{n=N+1}^{\infty}\left\Vert g_{n}\right\Vert
_{X}<\varepsilon/4.$ Then%
\[
\mathcal{D}^{\omega^{\alpha}}\left(  F,\varepsilon,X\right)  \subseteq
\mathcal{D}^{\omega^{\alpha}}\left(  \sum_{n=1}^{N}g_{n},\frac{\varepsilon}%
{2},X\right)  =\emptyset,
\]
since $\beta_{X}\left(  \sum_{n=1}^{N}g_{n}\right)  <\omega^{\alpha}$ by
Proposition \ref{P2}$.$ Thus $\beta_{X}\left(  F\right)  \leq\omega^{\alpha}.$
\end{proof}

\begin{corollary}
[{Kuratowski, \cite[\S 31, VI]{K}}]Let $X$ be a Polish space and $Y$ be a
$G_{\delta}$ subset of $X.$ Then every function of Baire class one on $Y$ can
be extended to a function of Baire class one on $X.$
\end{corollary}

\section{Extension of Continuous Functions}

In this section, we study the extension of a continuous function on a subspace
of a metric space to the whole space. To begin with, we consider the extension
of a\ continuous function from a dense subspace.

Consider a metric space $X$ with a dense subspace $Y$. Suppose that
$f:Y\rightarrow\mathbb{R}$ is continuous on $Y$. An obvious way of extending
$f$ to $X$ (if $f$ is loacally bounded) is to consider the limit superior (or
limit inferior) of $f,$ i.e.,
\[
\tilde{f}\left(  x\right)  =\limsup_{y\rightarrow x,\text{ }y\in Y}f\left(
y\right)  =\inf_{\delta>0}\sup_{\substack{d\left(  x,y\right)  <\delta\\y\in
Y}}f\left(  y\right)  .
\]
The extended function, which is upper semi-continuous (lower semi-continuous
in the case of limit inferior), is not optimal as far as the oscillation index
is concerned. In fact, the $\limsup$ extension $\tilde{f}$\ of the continuous
function $f$ in Example \ref{Ex1} below has oscillation index $\beta
_{X}\left(  \tilde{f}\right)  =\omega.$ The following is an alternative
algorithm that produces an extension with the smallest possible oscillation
index. If $A\subseteq\operatorname*{dom}f,$ $\operatorname*{osc}\left(
f,A\right)  $ is defined to be $\sup\left\{  \left\vert f\left(  x\right)
-f\left(  x^{\prime}\right)  \right\vert :x,x^{\prime}\in A\right\}  .$ If $x$
belongs to the closure of $\operatorname*{dom}f,$ then define%
\[
\operatorname*{osc}\left(  f,x\right)  =\lim_{\delta\rightarrow0}%
\operatorname*{osc}\left(  f,B\left(  x,\delta\right)  \cap\operatorname*{dom}%
f\right)  .
\]

We first define layers of approximate extensions inductively. Let $S_{0}=X$
and $n_{0}(s)=0$ for all $s\in S_{0}$.\ Assume that $S_{k}$ has been chosen
and $n_{k}(s)$ is defined for all $s\in S_{k}$. Let $\mathcal{U}%
_{k}=\{B(s,2^{-n_{k}(s)}):s\in S_{k}\}$ and $X_{k}=\cup\,\mathcal{U}_{k}$.
Choose a partition of unity $(\varphi_{s}^{k})_{s\in S_{k}}$ on $X_{k}$
subordinated to $\mathcal{U}_{k}$. For each $s\in S_{k}$, choose $y_{s}^{k}\in
Y\cap B(s,2^{-n_{k}(s)})$. Define $F_{k}:X_{k}\rightarrow\mathbb{R}$ by
$F_{k}(x)=\sum_{s\in S_{k}}\varphi_{s}^{k}(x)f(y_{s}^{k})$. For each $x\in
X_{k}$, let $S_{k}(x)=\{s\in S_{k}:x\in\operatorname{supp}\varphi_{s}^{k}\}$
and $l_{k}(x)=\max\{n_{k}(s):s\in S_{k}(x)\}+1$. Let $S_{k+1}$ be the set of
all $x\in X_{k}$ such that $\operatorname{osc}(f,x)<2^{-l_{k}(x)}$. If $x\in
S_{k+1}$, choose $n_{k+1}(x)\geq l_{k}(x)$ so that

\begin{enumerate}
\item $\operatorname{osc}(f,B(x,2^{1-n_{k+1}(x)})\cap Y)<2^{-l_{k}(x)}$,

\item \label{(2)}$B(x,2^{-n_{k+1}(x)})\subseteq B(s,2^{-n_{k}(s)})$ for all
$s\in S_{k}(x)$,

\item \label{(3)}$B(x,2^{1-n_{k+1}(x)})\cap\operatorname{supp}\varphi_{s}%
^{k}=\emptyset$ if $s\in S_{k}\backslash S_{k}(x)$.
\end{enumerate}

\noindent The extension $F$ (defined after Lemma \ref{L2.2}) is obtained by
pasting the layers $\left(  F_{k}\right)  $ one after another. Observe that
$Y\subseteq S_{k}\subseteq X_{k}$ for all $k$ and that \thinspace
$X_{k+1}\subseteq X_{k}$ because of condition (\ref{(2)}).

\begin{lemma}
\label{L2.1}Suppose that $s\in S_{k}$, $t\in S_{m}$ for some $m>k$, and that
$\operatorname{supp}\varphi_{s}^{k}\cap\operatorname{supp}\varphi_{t}^{m}%
\neq\emptyset$. Then $B(t,2^{-n_{m}(t)})\subseteq B(s,2^{-n_{k}(s)})$.
\end{lemma}

\begin{proof}
Let $x\in\operatorname{supp}\varphi_{s}^{k}\cap\operatorname{supp}\varphi
_{t}^{m}$. Then $x\in X_{j}$ for all $j\leq m$. In particular, if $m>j>k$,
then there exists $s_{j}\in S_{j}$ such that $x\in\operatorname{supp}%
\varphi_{s_{j}}^{j}$. Thus it suffices to prove the lemma for $m=k+1$. Assume
that $x\in\operatorname{supp}\varphi_{s}^{k}\cap\operatorname{supp}\varphi
_{t}^{k+1}$. Note that $s\in S_{k}(t)$. For otherwise, $B(t,2^{1-n_{k+1}%
(t)})\cap\operatorname{supp}\varphi_{s}^{k}=\emptyset$ by (\ref{(3)}). Since
$x$ belongs to this set, we have reached a contradiction. It now follows from
(\ref{(2)}) that $B(t,2^{-n_{k+1}(t)})\subseteq B(s,2^{-n_{k}(s)})$.
\end{proof}

\begin{lemma}
\label{L2.2}Suppose that $x\in X_{m}$ and $m>k\geq1$. Then there exists $s\in
S_{k}(x)$ such that $|F_{k}(x)-F_{m}(x)|<2^{1-l_{k-1}(s)}$. Moreover, if $x\in
Y$, then $|F_{k}(x)-f(x)|<2^{-l_{k-1}(s)}$ for some $s\in S_{k}(x)$.
\end{lemma}

\begin{proof}
Denote by $S$ the set of all $t\in S_{m}$ such that $\varphi_{t}^{m}(x)>0$ and
choose a point $y\in\cap_{t\in{S}}B(t,2^{-n_{m}(t)})\cap Y$. Let $s$ be an
element where $l_{k-1}(s)$ attains its minimum over $S_{k}\left(  x\right)  $.
By Lemma \ref{L2.1}, $B(t,2^{-n_{m}(t)})\subseteq B(s,2^{-n_{k}(s)})$ for all
$t\in S$. Hence $|f(y)-f(y_{t}^{m})|<2^{-l_{k-1}(s)}$ for any $t\in S$. By
Lemma \ref{L2.1} again, $y\in B(t,2^{-n_{m}(t)})\subseteq B(s^{\prime
},2^{-n_{k}(s^{\prime})})$ for all $t\in S$ and all $s^{\prime}\in S_{k}(x)$.
Hence
\[
|f(y)-f(y_{s^{\prime}}^{k})|<2^{-l_{k-1}(s^{\prime})}\leq2^{-l_{k-1}(s)}%
\]
for all $s^{\prime}\in S_{k}(x)$. Therefore
\begin{align*}
|F_{k}(x)-F_{m}(x)|  &  \leq|F_{k}(x)-f(y)|+|f(y)-F_{m}(x)|\\
&  <2^{-l_{k-1}(s)}+2^{-l_{k-1}(s)}=2^{1-l_{k-1}(s)}.
\end{align*}
Moreover, if $x\in Y$, then the above applies for $y=x$. Hence $|F_{k}%
(x)-f(x)|<2^{-l_{k-1}(s)}$.
\end{proof}

Observe that $l_{k}(s)\geq k+1$ for all $s\in S_{k}$, $k\geq0$. It follows
from Lemma \ref{L2.2} that $(F_{k})$ converges pointwise on $\cap X_{k}$ and
that the limit is $f$ on $Y$. Define $F:X\rightarrow\mathbb{R}$ by
\[
F(x)=%
\begin{cases}
\lim_{k}F_{k}(x) & \text{ if $x\in\cap X_{k}$},\\
F_{k}(x) & \text{ if $x\in X_{k}\backslash X_{k+1}$, $k\geq0$.}%
\end{cases}
\]
Then $F$ is an extension of $f$ to $X$.

\begin{lemma}
\label{L2.3}Suppose that $x\in X_{k}$ for some $k\geq1$. There exists an open
neighborhood $U$ of $x$ and $s\in S_{k}(x)$ such that
$|F(z)-F(x)|<2^{3-l_{k-1}(s)}$ for all $z\in U$.
\end{lemma}

\begin{proof}
Let $s$ be an element where $l_{k-1}(s)$ attains its minimum over
$S_{k}\left(  x\right)  $. Note that $F_{k}$ is continuous on the open set
$X_{k}$. Hence there is an open neighborhood $U$ of $x$ such that

\begin{enumerate}
\item $\operatorname{osc}(F_{k},U)<2^{-l_{k-1}(s)}$,

\item $U\subseteq X_{k}$,

\item $U\cap\operatorname{supp}\varphi_{s}^{k}=\emptyset$ if $s\in
S_{k}\backslash S_{k}(x)$.
\end{enumerate}

\noindent We claim that ${S}_{k}(z)\subseteq S_{k}(x)$ for all $z\in U$.
Indeed, if $z\in U$ and $s\in{S}_{k}(z)\backslash S_{k}(x)$, then $z\in
U\cap\operatorname{supp}\varphi_{s}^{k}=\emptyset$, a contradiction. Now if
$z\in U$, then either $z\in X_{m}$ for all $m$ or $z\in X_{m}\backslash
X_{m+1}$ for some $m\geq k$. In either case, $|F_{k}(z)-F(z)|\leq
2^{1-l_{k-1}(s)}$ by Lemma \ref{L2.2}. Therefore,
\begin{align*}
|F(z)-F(x)|  &  \leq|F(z)-F_{k}(z)|+|F_{k}(z)-F_{k}(x)|+|F_{k}(x)-F(x)|\\
&  <2^{1-l_{k-1}(s)}+2^{-l_{k-1}(s)}+2^{1-l_{k-1}(s)}<2^{3-l_{k-1}(s)}.
\end{align*}

\end{proof}

\noindent\ The next proposition is an immediate consequence of Lemma
\ref{L2.3}.

\begin{proposition}
Every $x\in\cap X_{k}$ is a point of continuity of $F$.
\end{proposition}

\begin{proposition}
If $x\in{\mathcal{D}}^{1}(F,2^{-m},X)\cap X_{k}$, $k\geq1$, then there exists
$s\in S_{k}(x)$ such that $l_{k-1}(s)\leq m+3$.
\end{proposition}

\begin{proof}
Since $x\in X_{k},$ by Lemma \ref{L2.3}, there exist an open neighborhood $U$
of $x$ and $s\in S_{k}\left(  x\right)  $ such that for all $z\in U,$
$\left\vert F\left(  z\right)  -F\left(  x\right)  \right\vert <2^{3-l_{k-1}%
\left(  s\right)  }.$ Hence $\left\vert F\left(  z_{1}\right)  -F\left(
z_{2}\right)  \right\vert <2^{4-l_{k-1}\left(  s\right)  }$ for all
$z_{1},z_{2}\in U.$ As $x\in{\mathcal{D}}^{1}(F,2^{-m},X),$ $-m<4-l_{k-1}%
\left(  s\right)  .$ Thus $l_{k-1}\left(  s\right)  \leq m+3.$
\end{proof}

\begin{proposition}
\label{P2.3}Suppose that $x\in X_{k}\cap{\mathcal{D}}^{2}(F,2^{-m},X)$,
$k\geq0$. Then $n_{k}(s)\leq m+2$ for all $s\in S_{k}$ such that $\varphi
_{s}^{k}(x)>0$.
\end{proposition}

\begin{proof}
Choose an open neighborhood $U_{1}$ of $x$ such that $U_{1}\subseteq
\{\varphi_{s}^{k}>0\}$ for all $s\in S_{k}$ such that $\varphi_{s}^{k}(x)>0$.
Note that, in particular, $U_{1}\subseteq X_{k}$. Then choose an open
neighborhood $U_{2}$ of $x$ such that $\operatorname{osc}(F_{k},U_{2})<2^{-m}%
$. Let $U=U_{1}\cap U_{2}$. There exist $z_{1},z_{2}\in U\cap{\mathcal{D}}%
^{1}(F,2^{-m},X)$ such that $|F_{k}(z_{1})-F_{k}(z_{2})|\geq2^{-m}$. If
$z_{1},z_{2}\notin X_{k+1}$, then $F(z_{i})=F_{k}(z_{i})$, $i=1,2$. This leads
to a contradiction with the fact that $\operatorname{osc}(F_{k},U_{2})<2^{-m}%
$. Thus at least one of $z_{1},z_{2}$ belongs to $X_{k+1}$. Denote it by $z$.
By the previous proposition, there exists $t\in S_{k+1}(z)$ such that
$l_{k}(t)\leq m+3$. Let $s\in S_{k}$ be such that $\varphi_{s}^{k}(x)>0$. We
claim that $s\in S_{k}(t)$. For otherwise, $B(t,2^{1-n_{k+1}(t)}%
)\cap\operatorname{supp}\varphi_{s}^{k}=\emptyset$. This is absurd since the
intersection contains the point $z$. It follows from that claim that
$l_{k}(t)\geq n_{k}(s)+1$. Hence $n_{k}(s)\leq m+2$, as required.
\end{proof}

\begin{proposition}
$\beta_{X}(F)\leq3$.
\end{proposition}

\begin{proof}
Suppose that $x\in{\mathcal{D}}^{3}(F,2^{-m},X)$ for some $m$. Then there
exists $k$ such that $x\in X_{k}\backslash X_{k+1}$. Choose a neighborhood $U$
of $x$ such that $U\subseteq B(x,2^{-m-2})\cap X_{k}$ and $\operatorname{osc}%
(F_{k},U)<2^{-m}$. There exist $z_{1},z_{2}\in U\cap{\mathcal{D}}^{2}%
(F,2^{-m},X)$ such that $|F(z_{1})-F(z_{2})|\geq2^{-m}$. If $z_{1},z_{2}\notin
X_{k+1}$, then $F(z_{i})=F_{k}(z_{i})$, $i=1,2$. This contradicts the fact
that $\operatorname{osc}(F_{k},U)<2^{-m}$. Hence there exists $z\in U\cap
X_{k+1}\cap{\mathcal{D}}^{2}(F,2^{-m},X)$. By Proposition \ref{P2.3},
$n_{k+1}(t)\leq m+2$ for all $t\in S_{k+1}$ such that $\varphi_{t}^{k+1}%
(z)>0$. Fix such a $t$. Note that
\begin{align*}
d(x,t)  &  \leq d(x,z)+d(z,t)\\
&  <2^{-m-2}+2^{-n_{k+1}(t)}\leq2^{1-n_{k+1}(t)}.
\end{align*}
Thus
\[
\operatorname{osc}(f,x)\leq\operatorname{osc}(f,B(t,2^{1-n_{k+1}(t)})\cap
Y)<2^{-l_{k}(t)}.
\]
We claim that $S_{k}(x)\subseteq S_{k}(t)$. For otherwise, there exists $s\in
S_{k}(x)\backslash S_{k}(t)$. Then $B(t,2^{1-n_{k+1}(t)})\cap
\operatorname{supp}\varphi_{s}^{k}=\emptyset$. This is absurd since the
intersection contains the point $x$. It follows from the claim that
$l_{k}(t)\geq l_{k}(x)$. Hence $\operatorname{osc}(f,x)<2^{-l_{k}(x)}$. Then
$x\in S_{k+1}\subseteq X_{k+1}$, a contradiction.
\end{proof}

We have shown that:

\begin{theorem}
\label{T2}Every continuous function $f$ on a dense subspace of a metric space
$X$ can be extended to a function $F$ on $X$ with $\beta_{X}\left(  F\right)
\leq3.$
\end{theorem}

\begin{theorem}
\emph{(\cite[Theorem 3.6]{LT})} Let $Y$ be a closed subspace of a metric space
$X$ and let $f$ be a function on $Y$ with $\beta_{Y}\left(  f\right)
<\omega_{1}.$ Then there exists a function $F$ on $X$ such that
\[
F_{|Y}=f\text{ and }\beta_{X}\left(  F\right)  =\beta_{Y}\left(  f\right)  .
\]

\end{theorem}

\begin{theorem}
\label{T3}Let $X$ be a metric space and $Y$ be a subspace of $X.$ Every
continuous function $f$ on $Y$ can be extended to a function $F$ on $X$ with
$\beta_{X}\left(  F\right)  \leq3.$
\end{theorem}

The following example shows that Theorem \ref{T3} is optimal.

\begin{example}
\label{Ex1}There is a subspace $Y\subseteq\left\{  0,1\right\}  ^{\omega}=X$
and a continuous real-valued function $f$ on $Y$ such that for any extension
$F$ of $f$ to $X,$ $\beta_{X}\left(  F\right)  \geq3.$
\end{example}

\begin{proof}
For any integer $n$, denote $n\left(  \operatorname{mod}2\right)  $ by
$\hat{n}.$ Let
\[
Y=\left\{  \left(  \varepsilon_{i}\right)  \in X:\varepsilon_{i}=0\text{ for
infinitely many }i\text{'s}\right\}  .
\]
We denote elements in $X$ of the form%
\[
\left(  \underset{n_{1}}{\underbrace{1,1,...,1}},0,\underset{n_{2}%
}{\underbrace{1,1,...,1}},0,...,\underset{n_{k}}{\underbrace{1,1,...,1}%
},0,...\right)
\]
by
\[
\left(  1^{n_{1}},0,1^{n_{2}},0,...,1^{n_{k}},0,...\right)  .
\]
Also write $\left(  \varepsilon_{1},...,\varepsilon_{k},\varepsilon
,\varepsilon,...\right)  $ as $\left(  \varepsilon_{1},...,\varepsilon
_{k},\varepsilon^{\omega}\right)  ,$ $\varepsilon_{i},\varepsilon\in\left\{
0,1\right\}  .$ Define $g:Y\rightarrow X$ by
\[
g\left(  1^{n_{1}},0,1^{n_{2}},0,...,1^{n_{k}},0,...\right)  =\left(  \hat
{n}_{1},\hat{n}_{2},...\right)  ,\text{ }n_{1},n_{2},...\in\mathbb{N\cup
}\left\{  0\right\}  ,
\]
and let $h:X\rightarrow\mathbb{R}$ be the canonical embedding of $X$ into
$\mathbb{R}$, $h\left(  \varepsilon_{1},\varepsilon_{2},...\right)
=\sum_{k=1}^{\infty}\frac{2\varepsilon_{k}}{3^{k}}.$ Then the function
$f=h\circ g:Y\rightarrow\mathbb{R}$ is continuous. Suppose that $F$ is an
extension of $f$ to $X$ such that $\beta_{X}\left(  F\right)  \leq2.$ First
observe that for any $n_{1},...,n_{k}\in\mathbb{N\cup}\left\{  0\right\}  $
and all $n\in\mathbb{N}$,%
\[
\left\vert F\left(  1^{n_{1}},0,...,1^{n_{k}},0,1^{2n},0^{\omega}\right)
-F\left(  1^{n_{1}},0,...,1^{n_{k}},0,1^{2n-1},0,1,0,1,...\right)  \right\vert
=\frac{1}{3^{k}}.
\]
Hence $\left(  1^{n_{1}},0,...,1^{n_{k}},0,1^{\omega}\right)  \in
\mathcal{D}^{1}\left(  F,\frac{1}{3^{k}},X\right)  .$ Let $F\left(  1^{\omega
}\right)  =a.$ Either $\left\vert a\right\vert \geq\frac{1}{2}$ or $\left\vert
1-a\right\vert \geq\frac{1}{2}.$ We assume the former; the proof is similar
for the latter case. Since $\left(  1^{\omega}\right)  \notin\mathcal{D}%
^{2}\left(  F,\frac{1}{3},X\right)  ,$ there exists a neighborhood $U$ of
$\left(  1^{\omega}\right)  $ such that $\left\vert F\left(  x\right)
-a\right\vert <\frac{1}{3}$ if $x\in U\cap\mathcal{D}^{1}\left(  F,\frac{1}%
{3},X\right)  .$ In particular, there exists $n_{1}\in\mathbb{N}$ such that
\[
\left\vert F\left(  1^{2n_{1}},0,1^{\omega}\right)  -a\right\vert =\frac{1}%
{3}-\delta\text{ for some }\delta>0.
\]
Similarly, using the fact that $\left(  1^{2n_{1}},0,1^{\omega}\right)
\notin\mathcal{D}^{2}\left(  F,\frac{1}{3^{2}},X\right)  ,$ we obtain
$n_{2}\in\mathbb{N}$ such that
\[
\left\vert F\left(  1^{2n_{1}},0,1^{2n_{2}},0,1^{\omega}\right)  -F\left(
1^{2n_{1}},0,1^{\omega}\right)  \right\vert <\frac{1}{3^{2}}.
\]
Continuing, we choose $n_{1},n_{2},...\in\mathbb{N}$ such that
\[
\left\vert F\left(  1^{2n_{1}},0,...,1^{2n_{k+1}},0,1^{\omega}\right)
-F\left(  1^{2n_{1}},0,...,1^{2n_{k}},0,1^{\omega}\right)  \right\vert
<\frac{1}{3^{k+1}},\text{ }k\in\mathbb{N}.
\]
In particular,%
\[
\left\vert F\left(  1^{2n_{1}},0,...,1^{2n_{k}},0,1^{\omega}\right)
-a\right\vert \leq\frac{1}{3}+\frac{1}{3^{2}}+...-\delta=\frac{1}{2}%
-\delta,\text{ }k\in\mathbb{N}.
\]
Since $\left\vert a\right\vert \geq\frac{1}{2},$ we have $\left\vert F\left(
1^{2n_{1}},0,...,1^{2n_{k}},0,1^{\omega}\right)  \right\vert \geq\delta$ for
all $k\in\mathbb{N}.$ But
\[
F\left(  1^{2n_{1}},0,...,1^{2n_{k}},0,1^{2n},0^{\omega}\right)  =f\left(
1^{2n_{1}},0,...,1^{2n_{k}},0,1^{2n},0^{\omega}\right)  =0
\]
for all $n\in\mathbb{N}.$ Hence $\left(  1^{2n_{1}},0,...,1^{2n_{k}%
},0,1^{\omega}\right)  \in\mathcal{D}^{1}\left(  F,\delta,X\right)  $ for all
$k\in\mathbb{N}.$ However, note that the sequence $\left(  \left(  1^{2n_{1}%
},0,...,1^{2n_{k}},0,1^{\omega}\right)  \right)  _{k\in\mathbb{N}}$ converges
to the point $\left(  1^{2n_{1}},0,...,1^{2n_{j}},0,1^{2n_{j+1}},0,...\right)
$ and
\begin{align*}
&  \left\vert F\left(  1^{2n_{1}},0,...,1^{2n_{k}},0,1^{\omega}\right)
-F\left(  1^{2n_{1}},0,...,1^{2n_{j}},0,1^{2n_{j+1}},0,...\right)  \right\vert
\\
&  =\left\vert F\left(  1^{2n_{1}},0,...,1^{2n_{k}},0,1^{\omega}\right)
-f\left(  1^{2n_{1}},0,...,1^{2n_{j}},0,1^{2n_{j+1}},0,...\right)  \right\vert
\\
&  =\left\vert F\left(  1^{2n_{1}},0,...,1^{2n_{k}},0,1^{\omega}\right)
\right\vert \geq\delta
\end{align*}
for all $n\in\mathbb{N}.$ Therefore, $\left(  1^{2n_{1}},0,...,1^{2n_{j}%
},0,1^{2n_{j+1}},0,...\right)  \in\mathcal{D}^{2}\left(  F,\delta,X\right)  ,$
contrary to the assumption that $\beta_{X}\left(  F\right)  \leq2.$
\end{proof}

Our final result presents a special class of spaces where the conclusion of
Theorem \ref{T3} may be improved upon. Recall that a topological space is
$0$-\emph{dimensional} if every open cover has a refinement that is an open
cover and consists of pairwise disjoint sets. In particular, a 0-dimensional
space has a basis consisting of clopen sets. Also note that a closed subspace
of a 0-dimensional space is 0-dimensional. If $A$ is a subset of a topological
space $X,$ the derived set $A^{\prime}$ of $A$ is the set of all cluster
points of $A.$ Let $A^{\left(  0\right)  }=A.$ If $A^{\left(  \alpha\right)
}$ has been defined, let $A^{\left(  \alpha+1\right)  }=\left(  A^{\left(
\alpha\right)  }\right)  ^{\prime}.$ If $\beta$ is a limit ordinal, let
\[
A^{\left(  \beta\right)  }=\bigcap\limits_{\alpha<\beta}A^{\left(
\alpha\right)  }.
\]
A topological space $X$ is said to be \emph{scattered} if $X^{\left(
\gamma\right)  }=\emptyset$ for some ordinal $\gamma.$

\begin{theorem}
Suppose that $Y$ is a subspace of a $0$-dimensional scattered metrizable space
$X$. If $f:Y\rightarrow\mathbb{R}$ is a continuous function, then there is an
extension $F:X\rightarrow\mathbb{R}$ of $f$ such that $\beta_{X}(F)\leq2$ and
that $F$ is continuous at every point in $Y$.
\end{theorem}

\begin{proof}
Since $X$ is scattered, $X^{(\gamma)}=\emptyset$ for some ordinal $\gamma.$
The proof is by induction on $\gamma.$ The case $\gamma=1$ is clear. Suppose
that the theorem holds for all $\gamma<\gamma_{0}$. Let $X$ be a
$0$-dimensional metrizable space with $X^{\left(  \gamma_{0}\right)
}=\emptyset$. For all $x\in X$, choose $\gamma_{x}<\gamma_{0}$ such that $x\in
X^{(\gamma_{x})}\backslash X^{(\gamma_{x}+1)}$. Let $d$ be a compatible metric
on $X$ that is bounded. Define $\delta_{x}=d(x,X^{(\gamma_{x})}\backslash
\{x\})$. Then $\delta_{x}>0$.\newline

\noindent\textbf{Case 1.} $\gamma_{0}$ is a limit ordinal.

\noindent Let $\mathcal{A}=\left\{  B\left(  x,\delta_{x}\right)  :x\in
X\right\}  .\,$Then $\mathcal{A}$ is an open cover of $X.$ Hence there is a
refinement $\mathcal{B}$ that is an open cover of $X$ consisting of pairwise
disjoint sets. In particular the elements of $\mathcal{B}$ are clopen subsets
of $X.$ If $U\in\mathcal{B}$, then $U\subseteq B\left(  x,\delta_{x}\right)  $
for some $x\in X.$ Hence $U\cap X^{\left(  \gamma_{x}+1\right)  }=\emptyset.$
Since $\gamma_{x}+1<\gamma_{0},$ we may apply the inductive hypothesis to
obtain an extension $f_{U}:U\rightarrow\mathbb{R}$ of $f_{|Y\cap U}$ such that
$\beta_{U}\left(  f_{U}\right)  \leq2$ and that $f_{U}$ is continuous at every
point in $Y\cap U.$ Take $F=\cup_{U\in\mathcal{B}}f_{U}.$ Since each $U$ is
clopen in $X,$ $F$ is continuous at each point in $Y$ and $\mathcal{D}%
^{2}\left(  F,\varepsilon,X\right)  \cap U=\mathcal{D}^{2}\left(
F,\varepsilon,U\right)  =\emptyset$ for all $\varepsilon>0$ and $U\in
\mathcal{B}.$ Therefore $\beta_{X}\left(  F\right)  \leq2.$\newline

\noindent\textbf{Case 2.} $\gamma_{0}$ is a successor ordinal.

For each $x\in X^{(\gamma_{0}-1)}$, choose a sequence $(W_{n,x})_{n=1}%
^{\infty}$ of clopen neighborhoods of $x$ such that $W_{n+1,x}\subseteq
W_{n,x}\subseteq B(x,1/n)$ for all $n\in\mathbb{N}$ and $W_{1,x}\subseteq
B(x,\delta_{x}/3)$. If $x$ and $x^{\prime}$ are distinct elements in
$X^{(\gamma_{0}-1)}$, then $B(x,\delta_{x}/3)\cap B(x^{\prime},\delta
_{x^{\prime}}/3)=\emptyset$. Hence $W_{1,x}\cap W_{1,x^{\prime}}=\emptyset$.
Note that $W_{1}=\cup\{W_{1,x}:x\in X^{(\gamma_{0}-1)}\}$ is clopen in $X$.
Indeed, clearly $W_{1}$ is open. If $z\in\overline{W_{1}}$, then choose
$(x_{n})$ in $X^{(\gamma_{0}-1)}$ and a sequence $(z_{n})$ converging to $z$
such that $z_{n}\in W_{1,x_{n}}$ for all $n$. If $(x_{n})$ has a constant
subsequence, then clearly $z\in W_{1}$. Otherwise, assume that all $x_{n}$'s
are distinct. For all distinct $n,m\in\mathbb{N}$,
\begin{align*}
\max(\delta_{x_{n}},\delta_{x_{m}})  &  \leq d(x_{n},x_{m})\leq d(x_{n}%
,z_{n})+d(z_{n},z_{m})+d(z_{m},x_{m})\\
&  <\delta_{x_{n}}/3+d(z_{n},z_{m})+\delta_{x_{m}}/3.
\end{align*}
Hence $\max(\delta_{x_{n}},\delta_{x_{m}})/3<d(z_{n},z_{m})$. Since $(z_{n})$
converges, $\delta_{x_{n}}\rightarrow0$. Then
\[
d(x_{n},z)\leq d(x_{n},z_{n})+d(z_{n},z)<\delta_{x_{n}}/3+d(z_{n}%
,z)\rightarrow0.
\]
Since the $x_{n}$'s are distinct elements in $X^{(\gamma_{0}-1)}$, $z\in
X^{(\gamma_{0})}$, contrary to the assumption. Hence $W_{1}$ is clopen in $X$.

Now $(X\backslash W_{1})^{(\gamma_{0}-1)}=\emptyset$. Hence by the inductive
hypothesis, there exists an extension $f_{0}:X\backslash W_{1}\rightarrow
\mathbb{R}$ of $f_{|Y\cap(X\backslash W_{1})}$ such that $\beta_{X\backslash
W_{1}}(f_{0})\leq2$ and that $f_{0}$ is continuous at every point in
$Y\cap(X\backslash W_{1})$.

For each $x\in X^{(\gamma_{0}-1)}$ and each $n\in\mathbb{N}$, set
$U_{n,x}=W_{n,x}\backslash W_{n+1,x}$. Then $U_{n,x}^{(\gamma_{0}%
-1)}=\emptyset$. By the inductive hypothesis, there exists an extension
$f_{n,x}:U_{n,x}\rightarrow\mathbb{R}$ of $f_{|Y\cap U_{n,x}}$ such that
$\beta_{U_{n,x}}(f_{n,x})\leq2$ and that $f_{n,x}$ is continuous at every
point in $Y\cap U_{n,x}$. Consider $y\in Y\cap U_{n,x}$. Choose a clopen
neighborhood $V_{y}$ of $y$ such that $V_{y}\subseteq U_{n,x}\cap
B(y,\min(\delta_{y}/3,1/n))$ and that
\[
|f_{n,x}(z)-f(y)|=|f_{n,x}(z)-f_{n,x}(y)|<\min(\delta_{y},1/n)
\]
for all $z\in V_{y}$. Set $V=\cup\{V_{y}:y\in Y\cap(W_{1}\backslash
X^{(\gamma_{0}-1)})\}$. If $x\in X^{(\gamma_{0}-1)}\backslash Y$, define
$F(x)=0$. If $y\in X^{(\gamma_{0}-1)}\cap Y$, define $F(y)=f(y)$. Then define
\[
F(z)=%
\begin{cases}
f_{0}(z) & \text{ if $z\notin W_{1}$}\\
f_{n,x}(z) & \text{ if $z\in V\cap U_{n,x}$ for some $x\in X^{(\gamma_{0}-1)}$
and $n\in\mathbb{N}$}\\
F(x) & \text{ if $z\in W_{1,x}\backslash V$ for some $x\in X^{(\gamma_{0}-1)}%
$}.
\end{cases}
\]
Since $X\backslash W_{1}$ and all $V_{y}$ are open in $X$, by the definition
of $F$, we see that $\mathcal{D}^{2}\left(  F,\varepsilon,X\right)  \cap
(V\cup(X\backslash W_{1}))=\emptyset$ and that $F$ is continuous at every
point in $Y\cap(V\cup(X\backslash W_{1}))$. Suppose $y\in Y\cap X^{(\gamma
_{0}-1)}$. If $y$ is not a point of continuity of $F$, then there exists a
sequence $(z_{m})$ converging to $y$ such that $(F(z_{m}))$ does not converge
to $f(y)$. Without loss of generality, we may assume that $z_{m}\in W_{1,y}$
for all $m$. Since $F=f(y)$ on $W_{1,y}\backslash V$, we may also assume that
$z_{m}\in V$ for all
$m$. Choose sequences $(n_{m})$ in $\mathbb{N}$ and $(y_{m})$ in $Y$ so that
$y_{m}\in U_{n_{m},y}$ and $z_{m}\in V_{y_{m}}$ for all $m$. Since $(z_{m})$
converges to $y$, $(n_{m})$ diverges to $\infty$. Then $d(z_{m},y_{m}%
)<\min(\delta_{y_{m}}/3,1/n_{m})\rightarrow0$ and
\[
|F(z_{m})-f(y_{m})|=|f_{n_{m},y}(z_{m})-f(y_{m})|<\min(\delta_{y_{m}}%
,1/n_{m})\rightarrow0.
\]
Hence $(y_{m})$ converges to $y$ and $(f(y_{m}))$ converges to $f(y)$ since
$f$ is continuous on $Y$. But then $(F(z_{m}))$ converges to $f(y)$, a
contradiction. Hence $F$ is continuous at every point in $Y\cap X^{(\gamma
_{0}-1)}$ as well. Since $Y\subseteq V\cup(X\backslash W_{1})\cup
X^{(\gamma_{0}-1)}$, $F$ is continuous at all points in $Y$.

Finally, suppose that $z\in\mathcal{D}^{2}\left(  F,\varepsilon,X\right)  $
for some $\varepsilon>0$. By the above, $z\in W_{1}\backslash V$. Choose $x\in
X^{(\gamma_{0}-1)}$ such that $z\in W_{1,x}$. Then $F(z)=F(x)$. Choose
$(z_{m})$ in $W_{1,x}\cap\mathcal{D}^{1}\left(  F,\varepsilon,X\right)  $ such
that $(z_{m})$ converges to $z$ and $|F(z_{m})-F(z)|\geq\varepsilon/2$ for all
$m$. In particular, $z_{m}\in V$ for all $m$. Choose sequences $(n_{m})$ in
$\mathbb{N}$ and $(y_{m})$ in $Y$ so that $y_{m}\in U_{n_{m},x}$ and $z_{m}\in
V_{y_{m}}$ for all $m$. \newline

\noindent\textbf{Claim}. $\delta_{y_{m}}\rightarrow0$.\newline

\noindent If the claim fails, then by going to a subsequence if necessary, we
may assume that there exists $\delta>0$ such that $\delta_{y_{m}}\geq\delta$
for all $m$, that $d(z_{m},z_{k})<\delta/6$ for all $m,k$ and that
$(\gamma_{y_{m}})$ is a nondecreasing sequence of ordinals. If $m<k$ and
$y_{m}$ and $y_{k}$ are distinct, then $y_{k}\in X^{(\gamma_{y_{m}}%
)}\backslash\{y_{m}\}$. Thus
\begin{align*}
\delta_{y_{m}}  &  \leq d(y_{m},y_{k})\\
&  \leq d(y_{m},z_{m})+d(z_{m},z_{k})+d(z_{k},y_{k})\\
&  <\delta_{y_{m}}/3+\delta/6+\delta_{y_{k}}/3\\
&  \leq\delta_{y_{m}}/2+\delta_{y_{k}}/3.
\end{align*}
Hence $\delta_{y_{k}}>3\delta_{y_{m}}/2$ whenever $k>m$ and $y_{k}\neq y_{m}$.
Since the metric $d$ is assumed to be bounded, the sequence $(y_{m})$ must
have a constant subsequence. Without loss of generality, let $y_{m}=y_{0}$ for
all $m$. Then $z_{m}\in V_{y_{0}}$ for all $m$ and hence $z\in V_{y_{0}%
}\subseteq V$, a contradiction. This proves the claim.\newline

Using the claim, choose $m$ large enough so that $\delta_{y_{m}}%
<\varepsilon/2$. Now
\[
|F(v)-f(y_{m})|<\delta_{y_{m}}<\varepsilon/2
\]
for all $v\in V_{y_{m}}$. Since $V_{y_{m}}$ is a neighborhood of $z_{m}$, we
see that $z_{m}\notin\mathcal{D}^{1}\left(  F,\varepsilon,X\right)  $,
contrary to the choice of $z_{m}$.
\end{proof}


\begin{thebibliography}{9}                                                                                                %


\bibitem {B}\textsc{R. L. Blair,} Extensions of continuous functions from
dense subspaces. Proc. Amer. Math. Soc. \textbf{54} (1976), 355--359.

\bibitem {BBT}\textsc{A. M. Bruckner, J. B. Bruckner and B. S. Thomson,} Real
Analysis, Prentice-Hall, Inc., New Jersey, 1997.

\bibitem {D}\textsc{J. Dugundji}, Topology, Allyn and Bacon, Inc., Boston, 1966.

\bibitem {H}\textsc{S. Hern\'{a}ndez}, Extension of continuous functions into
uniform spaces. Proc. Amer. Math. Soc. \textbf{97} (1986), no. 2, 355--360.

\bibitem {KL}\textsc{A. S. Kechris and A. Louveau}, A classification of
Baire-1 functions, Trans. Amer. Math. Soc. \textbf{318}(1990), 209--236.

\bibitem {KS}\textsc{O. Kalenda and J. Spurny, }Extending Baire-one functions
on topological spaces, Topology and its Appl., \emph{to appear.}

\bibitem {K}\textsc{C. Kuratowski,} Topologie, Vol. I. (French) 4\`{e}me
\'{e}d. Monografie Matematyczne, Tom 20. Pa\'{n}stwowe Wydawnictwo Naukowe,
Warsaw, 1958.

\bibitem {LT}\textsc{D. H. Leung and W.-K. Tang,} Functions of Baire Class
One, Fund. Math., \textbf{179}(2003), 225--247.
\end{thebibliography}
\end{document}